\documentclass{article}

\usepackage{bera}
\usepackage[T1]{fontenc}
\RequirePackage{geometry}
\usepackage[colorlinks,citecolor=blue, linkcolor = blue,urlcolor=blue]{hyperref}
\usepackage{amsthm,amsmath,amsfonts,amssymb,mathtools}
\usepackage[numbers]{natbib}
\usepackage{authblk} 
\usepackage{graphicx}
\usepackage{enumerate}
\usepackage{tikz}
\usetikzlibrary{calc,arrows}
\usetikzlibrary{decorations.pathreplacing}
\usepackage{standalone}
\usepackage{bbm}

\newtheorem{theorem}{Theorem}[section]
\newtheorem{proposition}[theorem]{Proposition}
\newtheorem{lemma}[theorem]{Lemma}
\newtheorem{corollary}[theorem]{Corollary}
\newtheorem{remark}{Remark}[section]

\newtheorem{example}{Example}[section]

\newcommand\extrafootertext[1]{%
    \bgroup
    \renewcommand\thefootnote{\fnsymbol{footnote}}%
    \renewcommand\thempfootnote{\fnsymbol{mpfootnote}}%
    \footnotetext[0]{#1}%
    \egroup
}

\DeclareRobustCommand{\stirling}{\genfrac\{\}{0pt}{}}
\newcommand{\interior}[1]{%
  {\kern0pt#1}^{\mathrm{o}}%
}

\title{A note on maximal conditional entropy on Lebesgue spaces}

\author{Michael Hediger\footnote{University of Zurich, Switzerland. E-mail: \texttt{michael.hediger@math.uzh.ch}}}

\affil{Department of Mathematics, Winterthurerstrasse 190, 8057 Zurich, Switzerland}

\begin{document}

\maketitle

\begin{abstract}
Let $(X,\mathcal{B},P)$ be a probability space and $\mathit{a}$ be a sub $\sigma$-field that is generated by an increasing sequence of sub $\sigma$-fields $(\mathit{a}_{n})_{n \in \mathbb{N}}$. Given $\theta \in \Theta$, where $\Theta$ is some set, let $(X_{n}^{\theta})_{n \in \mathbb{N}}$ be a martingale adapted to $(\mathit{a}_{n})_{n \in \mathbb{N}}$. Martin (1969) provides sufficient conditions to show that $(X_{n}^{\theta})_{n \in \mathbb{N}}$ converges a.s.\ uniformly on $\Theta$ to a random variable $X^{\theta}$. His results are based on the assumption that there exists an integer $n$ s.t.\ the conditional entropy given $\mathit{a}_{n}$ is uniformly bounded over the set of finite partitions of $X$ with atoms from $\mathit{a}$. This study complements Martin's results by studying the latter assumption on the maximal conditional entropy in the context of measurable partitions of Lebesgue spaces. We provide conditions under which $\mathit{a}$ conveys too much information for the maximal conditional entropy to be finite. As an example, we consider the space of continuous functions with a compact support, equipped with the Borel $\sigma$-field.\\[2mm]
\end{abstract}

\section{Primary notation} \label{se:notation}
The set of strictly positive integers is denoted with $\mathbb{N} = \{n \in \mathbb{Z} \colon n \geq 1\}$. We write $(a_{n}) = (a_{n})_{n \in \mathbb{N}}$ for a sequence indexed over $\mathbb{N}$. Given a subset $A$ of a topological space $(X, \tau)$, the sets $\interior{A}$ and $\overline{A}$ denote the topological interior and closure of $A$, respectively. The Borel $\sigma$-field on $(X, \tau)$ is denoted with $\mathfrak{B}(X)$. If $x$ is an element of a metric space, $B(x,r)$ (resp.\ $B[x,r]$) denotes the open (resp.\ closed) ball with radius $r$ centered around $x$. Given two sets $A$ and $B$ of a measure space, we write $A \, \triangle \, B = A \setminus B \cup B \setminus A$ for the symmetric difference between $A$ and $B$. If $\mu_{1}$ and $\mu_{2}$ are two measures, $\mu_{2} \ll \mu_{1}$ means that $\mu_{2}$ is absolutely continuous with respect to $\mu_{1}$. Orthogonal measures are denoted with $\mu_{1} \perp \mu_{2}$. In the following we assume that any given measure space is complete in the measure theoretic sense.  

\section{Introduction} \label{se:intro}

Let $(X,\mathcal{B},P)$ be a probability space and $\mathit{a} \subset \mathcal{B}$ be a sub $\sigma$-field of $\mathcal{B}$. It is known that if $\mathit{a} = \sigma(\cup_{n \in \mathbb{N}}\mathit{a}_{n})$ for some increasing sequence of sub $\sigma$-fields $(\mathit{a}_{n})$, then the sets from $\mathit{a}$ can be approximated pointwise with sets from $\cup_{n \in \mathbb{N}}\mathit{a}_{n}$. That is, for any $\varepsilon > 0$ and any set $A \in \mathit{a}$, there exists $B \in \cup_{n \in \mathbb{N}}\mathit{a}_{n}$ s.t.\ $P(A \, \triangle \, B) < \varepsilon$. Let $Z_{0}^{\mathit{a}}$ be the sets of finite partitions of $X$ with atoms from $\mathit{a}$. In order to approximate sets from $\mathit{a}$ with sets from $\cup_{n \in \mathbb{N}}\mathit{a}_{n}$ uniformly over $\mathbb{N}$, Martin (1969) \cite{martin} introduces the assumption that there exists $n \in \mathbb{N}$ s.t.\ the conditional entropy given $\mathit{a}_{n}$ is uniformly bounded on $Z_{0}^{\mathit{a}}$.
\begin{lemma}[Lemma~1 in \cite{martin}] \label{lemma:martin}
Let $(\mathit{a}_{n})$ be an increasing sequence of sub $\sigma$-fields s.t.\ $\mathit{a} = \sigma(\cup_{n \in \mathbb{N}}\mathit{a}_{n})$. If there exists $n \in \mathbb{N}$ s.t.\ 
\begin{equation} \label{eq:martin}
    \sup_{\xi \in Z_{0}^{\mathit{a}}}E\bigg[-\sum_{A \in \xi}P(A\mid \mathit{a}_{n})\log P(A\mid \mathit{a}_{n}) \bigg] < \infty,
\end{equation}
then for any $\varepsilon > 0$ there exists $N \in \mathbb{N}$ s.t.\ for any $A \in \mathit{a}$ there exists $B \in \mathit{a}_{N}$ s.t.\ $P(A \, \triangle \, B) < \varepsilon$. 
\end{lemma}
The aim of the present article is to provide a better understanding of condition \eqref{eq:martin}. In particular, we study \eqref{eq:martin} in the context of measurable partitions of Lebesgue spaces (Section~\ref{se:recap} gives an introduction of terms). If the original space $(X,\mathcal{B},P)$ is a Lebesgue space and $\mathit{a} = \mathcal{B}$, we give necessary and sufficient conditions for \eqref{eq:martin} to be satisfied (Theorem~\ref{thm1}). If $\mathit{a}$ is not necessarily equal to $\mathcal{B}$, we provide conditions that enable us to deduce when \eqref{eq:martin} is not satisfied (Theorem~\ref{thm2}). Our arguments are accompanied with several examples. Specifically, we consider the case where $X = C(T)$, the space of real valued and continuous functions defined on a compact subset $T \subset \mathbb{R}^{d}$ and $\mathcal{B} = \mathfrak{B}(C(T))$, the Borel $\sigma$-field on $C(T)$. As shown in \cite{martin}, \eqref{eq:martin} is linked with the uniform convergence of martingale families (related results are given in \cite{boylan} and \cite{blake}). We add to the results of \cite{martin} and provide an extension in Section~\ref{se:martingales} (Corollary~\ref{corollary3}). An interesting example where \eqref{eq:martin} is not met is given in Section~\ref{se:UniformConvergence2}. There, the martingale family is obtained from the Radon–Nikodym derivative of the finite dimensional distributions of a given random field.

\section{Measurable partitions and Lebesgue spaces} \label{se:recap}
The following is based on Sections~1.2 and 1.3 of \cite{entropy}. A partition  $\xi$ of $(X,\mathcal{B},P)$ is a collection of subsets of $X$ which are disjoint and whose union is $X$. The elements of a partition $\xi$ are called atoms. As an example, $\epsilon = \{\{x\} \colon x \in X\}$ is referred to as the point partition on $(X,\mathcal{B},P)$. If $\xi$ is a partition of $(X,\mathcal{B},P)$, any subset of $X$ which is a union of atoms of $\xi$ is called an $\xi$-set.  A partition $\xi$ of $(X,\mathcal{B},P)$ is called measurable if there exists a countable collection $\mathcal{Q} = \{Q_{j} \colon j \in \mathbb{N}\}$ of measurable $\xi$-sets that separate the points of $\xi$, i.e., for any pair of atoms $A_{1},A_{2} \in \xi$ ($A_{1} \neq A_{2}$), there exists $Q \in \mathcal{Q}$ s.t.\ either $A_{1} \subset Q$ and $A_{2} \subset X\setminus Q$ or vice versa. The separating set $\mathcal{Q}$ associated with a measurable partiation $\xi$ of $(X,\mathcal{B},P)$ is referred to as a basis for $\xi$ (cf.\ Section~1.3 of \cite{rokhlin}). Let $J = \mathbb{N}$ or $J = \{1, \dotsc, N\}$, $N \in \mathbb{N}$. A countable collection $\Gamma = \{G_{j} \colon j \in J\}$ of measurable sets is called a basis for $(X,\mathcal{B},P)$ if $\sigma(\Gamma) = \mathcal{B}$ and $\Gamma$ separates the points of $X$, i.e., for any $x,y \in X$ ($x \neq y$) there exists $G \in \Gamma$ s.t.\ either $x \in G$ and $y \in X \setminus G$ or vice versa. Given $\alpha \in \{0,1\}^{J}$, let 
\begin{gather*}
    G_{j}^{\alpha(j)} = \begin{cases} G_{j}, & \alpha(j) = 0, \\ X \setminus G_{j}, & \alpha(j) = 1. \end{cases}
\end{gather*}
The probability space $(X,\mathcal{B},P)$ is said to be a Lebesgue space if it has a basis $\Gamma$ and if the partition built upon the atoms
\begin{gather*}
    A_{\alpha} = \bigcap_{j \in J}G_{j}^{\alpha(j)}, \quad  \alpha \in \{0,1\}^{J},
\end{gather*}
is equal to $\epsilon$, the point partition on $(X,\mathcal{B},P)$. Several examples of Lebesgue spaces are discussed in Section~1.2 of \cite{entropy}. 

For the rest of this section we shall assume that our initial space $(X,\mathcal{B},P)$ is in fact a Lebesgue space. Let $Z$ denote the the collection of all measurable partitions of $(X,\mathcal{B},P)$. Given $\xi, \eta \in Z$, we say that $\eta$ refines $\xi$ ($\xi \leq \eta$) if each atom of $\xi$ is a $\eta$-set. If $\xi \leq \eta$ and $\eta \leq \xi$ we say that $\xi$ and $\eta$ are equal and write $\xi = \eta$. Notice that the statement $\xi \leq \eta$ is understood up to sets of measure zero (compare to p.\ 10 of \cite{entropy}). It is known that the pair $(Z,\leq)$ is a complete lattice with infimum $\wedge_{i \in I}\xi_{i}$ and supremum $\vee_{i \in I}\xi_{i}$, where $I$ is some set and $\{\xi_{i} \colon i \in I\} \subset Z$. In particular, $\xi_{n} \uparrow \xi$ means that $\xi_{n} \leq \xi_{n+1}$ and $\xi = \vee_{n \in \mathbb{N}}\xi_{n}$. The case $\xi_{n} \downarrow \xi$ is treated analogously. Actually, any element of $Z$ can be written as the limit of an increasing sequence of finite and measurable partitions of $(X,\mathcal{B},P)$. Explicitly, let $Z_{0} \coloneqq \{\xi \in Z \colon \text{$\xi$ is finite}\}$. Then, for any $\xi \in Z$ there exists $(\xi_{n})$ s.t.\ $\xi_{n} \in Z_{0}$ for any $n \in \mathbb{N}$ and $\xi_{n} \uparrow \xi$ (see for instance Section~1.3 of \cite{rokhlin}). Notice that in the latter statement, the atoms of $\xi_{n}$ are built upon finite intersections of elements from a basis of $\xi$.

Finally, we state the definition of the conditional entropy of $\xi$ given $\eta$, where $\xi, \eta \in Z$. To do so, recall that since $(X,\mathcal{B},P)$ is a Lebesgue space, there exists a one to one correspondence $\varsigma$ between elements of $Z$ and sub $\sigma$-fields of $\mathcal{B}$. Given $\xi \in Z$, we let $\varsigma(\xi)$ denote the sub $\sigma$-field associated with $\xi$. The conditional entropy $H(\xi \mid \eta)$ of $\xi$ given $\eta$ is defined by 
\begin{gather*}
    H(\xi \mid \eta) = \begin{cases} E\bigg[-\sum_{A \in \xi}P\big(A\mid \varsigma(\eta)\big)\log P\big(A\mid \varsigma(\eta)\big) \bigg], & \text{if $\xi$ is countable,} \\ \infty, & \text{otherwise}. \end{cases}
\end{gather*}

\section{Maximal conditional entropy on Lebesgue spaces} \label{se:MaxCondEntropy}
If $(X,\mathcal{B},P)$ is a Lebesgue space, \eqref{eq:martin} of Lemma~\ref{lemma:martin} reads as follows: There exists $n \in \mathbb{N}$ s.t.\
\begin{gather*}
    \sup_{\xi \in Z_{0}^{\mathit{a}}}H\big(\xi \mid \varsigma^{-1}(\mathit{a}_{n})\big) < \infty,
\end{gather*}
where the elements of $Z_{0}^{\mathit{a}}$ are measurable partitions of $(X,\mathcal{B},P)$ since they are finite and their atoms are members of $\mathit{a}$. In that case $Z_{0}^{\mathit{a}} = \{\xi \in Z_{0} \colon \text{$\xi$ has atoms from $\mathit{a}$}\}$ and $Z_{0} = Z_{0}^{\mathcal{B}}$.

\begin{theorem} \label{thm1}
    Suppose that $(X,\mathcal{B},P)$ is a Lebesgue space and let $\eta$ be any measurable partition of $(X,\mathcal{B},P)$, i.e., $\eta \in Z$. Then, if
    \begin{enumerate}[(i)]
        \item \label{item:t1a} $X$ is finite, $\sup_{\xi \in Z_{0}} H( \xi \mid \eta) < \infty$;
        \item \label{item:t1b} $X$ is countably infinite, $\sup_{\xi \in Z_{0}} H( \xi \mid \eta) < \infty$ if and only if $H( \epsilon \mid \eta) < \infty$;
        \item \label{item:t1c} $X$ is not countable, $\sup_{\xi \in Z_{0}} H( \xi \mid \eta) = \infty$.
    \end{enumerate}
\end{theorem}

\begin{proof}
    The first item is clear since for any $\xi \in Z_{0}$, $H(\xi \mid \eta) < \infty$ and as $X$ is finite, the number of partitions of $X$ is finite as well and given by the bell number $\sum_{k=0}^{N}\stirling{N}{k}$, where $N$ is the number of elements in $X$ and $\stirling{N}{k}$ denotes the Stirling number of the second kind. To prove (\ref{item:t1b}), we recall that since $(X,\mathcal{B},P)$ is a Lebesgue space, the finest and measurable partition is given by the point partition $\epsilon$ on $(X,\mathcal{B},P)$. That is, for any $\xi \in Z$, $\xi \leq \epsilon$. Further, since $X$ is countable, $\epsilon$ is a countable as well. Therefore, using Corollary 2.6 of \cite{entropy}, for any $\xi \in Z_{0}$, $H(\xi \mid \eta) \leq H(\epsilon \mid \eta)$. Hence, if $H(\epsilon \mid \eta) < \infty$, then also $\sup_{\xi \in Z_{0}} H( \xi \mid \eta) < \infty$. To prove the other direction, assume by contradiction that $\sup_{\xi \in Z_{0}} H( \xi \mid \eta) < \infty$ but $H(\epsilon \mid \eta) = \infty$. Then, since $\epsilon$ is measurable, there exists $(\epsilon_{n})$ s.t.\ $\epsilon_{n} \in Z_{0}$ for any $n \in \mathbb{N}$ and $\epsilon_{n} \uparrow \epsilon$. By Theorem~2.19 of \cite{entropy}, $H(\epsilon_{n} \mid \eta) \uparrow H(\epsilon \mid \eta)$. But then,
    \begin{gather*}
        H(\epsilon \mid \eta) = \sup_{n \in \mathbb{N}}H(\epsilon_{n} \mid \eta) \leq \sup_{\xi \in Z_{0}} H( \xi \mid \eta) < \infty,
    \end{gather*}
    which gives the desired contradiction. Finally, let us show (\ref{item:t1c}). First, since $X$ is not countable, $\epsilon$ is a measurable partition s.t.\ $H(\epsilon \mid \eta) = \infty$. Suppose by contradiction that $\sup_{\xi \in Z_{0}} H( \xi \mid \eta) < \infty$. Again, we choose $(\epsilon_{n})$ s.t.\ $\epsilon_{n} \in Z_{0}$ for any $n \in \mathbb{N}$ and $\epsilon_{n} \uparrow \epsilon$. Then,
    \begin{gather*}
        \infty = H(\epsilon \mid \eta) = \sup_{n \in \mathbb{N}}H(\epsilon_{n} \mid \eta) \leq \sup_{\xi \in Z_{0}} H( \xi \mid \eta),
    \end{gather*}
    which gives a contradiction. This shows (\ref{item:t1c}) and concludes the proof of the theorem.
\end{proof}
Having in mind condition \eqref{eq:martin} of Lemma~\ref{lemma:martin}, let us assume that $\mathit{a} = \mathcal{B}$ and $(X,\mathcal{B},P)$ is a Lebesgue space. Then, Theorem~\ref{thm1} shows that if $X$ is not countable, \eqref{eq:martin} can not be satisfied. If $X$ is countably infinite, \eqref{eq:martin} is satisfied if and only if there exists $n \in \mathbb{N}$ s.t.\ $H( \epsilon \mid \varsigma^{-1}(\mathit{a}_{n})) < \infty$. 

\begin{example} \label{example1}
Let $C(T) = \{f \colon f \colon T \to \mathbb{R}, \text{ $f$ continuous}\}$, where $T \subset \mathbb{R}^{d}$, compact. We view $C(T)$ as a metric space, with metric $d(f,g) = \sup_{v \in T}\lvert f(v) - g(v) \rvert$ and equip $C(T)$ with the Borel $\sigma$-field $\mathfrak{B}(C(T))$ and a probability measure $P$. Then, the probability space $(C(T), \mathfrak{B}(C(T)), P)$ is a Lebesgue space (cf.\ Example 1.9 in \cite{entropy}). Given $n \in \mathbb{N}$, let 
\begin{gather*}
    \pi_{v_{1}, \dotsc, v_{n}} \colon C(T) \to \mathbb{R}^{n}, \quad v_{1}, \dotsc, v_{n} \in T,
\end{gather*}
be the coordinate map on $C(T)$, i.e., $\pi_{v_{1}, \dotsc, v_{n}}(f) = (f(v_{1}), \dotsc, f(v_{n}))$. Consider a collection of points $S = \{s_{i} \colon s_{i} \in T, i \in \mathbb{N}\}$ s.t.\ $\overline{S} = T$. The following lemma is known (see for instance \cite{StochasticProcesses}, Chapter~V). 
\begin{lemma} \label{lemma1}
    $\mathfrak{B}(C(T)) = \sigma(\cup_{n \in \mathbb{N}}\mathit{a}_{n})$, where 
    \begin{gather*}
        \mathit{a}_{n} = \sigma\big(\{\pi_{s_{1}, \dotsc, s_{n}}^{-1}(B_{n}) \colon B_{n} \in \mathfrak{B}(\mathbb{R}^{n})\}\big).
    \end{gather*}
\end{lemma}
Notice that $\mathit{a}_{n} \subset \mathit{a}_{n+1}$ and thus, upon the previous lemma we are in the setting of Lemma~\ref{lemma:martin}. Clearly $C(T)$ is not countable and hence we apply Theorem~\ref{thm1} to conclude that there does not exist $n \in \mathbb{N}$ s.t.\ $\sup_{\xi \in Z_{0}} H( \xi \mid \varsigma^{-1}(\mathit{a}_{n})) < \infty$, i.e., \eqref{eq:martin} is not satisfied. 
\end{example}

In the previous example, we were in the case where $\mathit{a} = \mathcal{B}$. We readily see that if $\mathit{a}$ is strictly smaller than $\mathcal{B}$, (\ref{item:t1a}) of Theorem~\ref{thm1} is still true with $Z_{0}$ replaced with $Z_{0}^{\mathit{a}}$. This is not the case for item (\ref{item:t1c}). As a trivial example, we can take $X$ not countable and $\mathit{a} = \{X,\emptyset\}$, the trivial $\sigma$-field on $X$. The next theorem facilitates a simple negation of \eqref{eq:martin} when $\mathit{a}$ is not necessarily equal to $\mathcal{B}$.

\begin{theorem} \label{thm2}
    Let $(X,\mathcal{B},P)$ be a Lebesgue space and $\eta \in Z$. Suppose that there exists an uncountable subspace $E \subset X$ s.t.\ $E \in \mathit{a}$ and $\Gamma_{E} = \{E_{j} \colon j \in \mathbb{N}\}$ is a basis for $(E, \mathit{a}, P_{\mathit{a}})$, where $P_{\mathit{a}}$ is a probability measure on $\mathit{a}$. Then, $\sup_{\xi \in Z_{0}^{\mathit{a}}} H( \xi \mid \eta) = \infty$.
\end{theorem}

\begin{proof}
    Let $\epsilon_{E}$ be the partition on $E$ with atoms, 
    \begin{gather*}
       A_{\alpha} = \bigcap_{j \in \mathbb{N}}E_{j}^{\alpha(j)}, \quad  \alpha \in \{0,1\}^{\mathbb{N}}, \quad 
    \end{gather*}
    where 
    \begin{gather*}
       E_{j}^{\alpha(j)} = \begin{cases} E_{j}, & \alpha(j) = 0, \\ E \setminus E_{j}, & \alpha(j) = 1. \end{cases}
    \end{gather*}
    Since $\Gamma_{E}$ is a basis for $(E, \mathit{a}, P_{\mathit{a}})$, it separates the points of $E$. In particular, $\epsilon_{E}$ is a measurable partition of $(E, \mathit{a}, P_{\mathit{a}})$ which has at least as many atoms as there are elements in $E$. Let
    \begin{gather*}
        \epsilon^{*} = \big\{\{X \setminus E\} \cup \epsilon_{E}\big\}.
    \end{gather*}
    Then, $\epsilon^{*}$ is a measurable partition of $(X,\mathcal{B},P)$. To see it, we let $\mathcal{Q}_{E} = \{Q_{j} \colon j \in \mathbb{N}\}$ be a basis for $\epsilon_{E}$ and define the collection
    \begin{gather*}
        \mathcal{Q}^{*} = \big\{\{X \setminus E\} \cup \mathcal{Q}_{E} \big\}.
    \end{gather*}
    Notice that the members of $\mathcal{Q}^{*}$ are $\epsilon^{*}$-sets, members of $\mathcal{B}$ and $\mathcal{Q}^{*}$ separates the atoms of $\epsilon^{*}$. Thus, since $\mathcal{Q}^{*}$ is countable, it is a basis for $\epsilon^{*}$ and hence $\epsilon^{*} \in Z$. Therefore, since $\epsilon^{*}$ is a measurable partition of a Lebesgue space, there exists $(\epsilon^{*}_{n})$ s.t.\ $\epsilon^{*}_{n} \in Z_{0}$ for any $n \in \mathbb{N}$ and $\epsilon^{*}_{n} \uparrow \epsilon^{*}$. Recall also, that for a given $n \in \mathbb{N}$, the atoms of $\epsilon^{*}_{n}$ are built upon finite intersections of sets from $\mathcal{Q}^{*}$, which are all members of $\mathit{a}$. Thus, $\epsilon^{*}_{n} \in Z_{0}^{\mathit{a}}$, $n \in \mathbb{N}$. In conclusion, since $E$ is not countable, $\sup_{n \in \mathbb{N}}H( \epsilon^{*}_{n} \mid \eta) = H(\epsilon^{*} \mid \eta) = \infty$, and the proof is complete.
\end{proof}

\begin{example} \label{example2}
We remain in the setting of Example~\ref{example1} and consider the Lebesgue space $(C(T), \mathfrak{B}(C(T)), P)$. Let $F$ be a closed subset of $T$ with positive Lebesgue measure. Consider a collection of points $S = \{s_{i} \colon s_{i} \in \interior{F}, i \in \mathbb{N}\}$ s.t.\ $\overline{S} = F$. We fix $h \in C(T)$ and consider the sub $\sigma$-field $\mathit{a} = \sigma(\{G_{F} \cup \{E\}\})$, where 
\begin{gather*}
    G_{F} = \{\pi_{s_{1}, \dotsc, s_{n}}^{-1}(B_{n}) \colon B_{n} \in \mathfrak{B}(\mathbb{R}^{n}), n \in \mathbb{N}\},
\end{gather*}
and $E = \{f \in C(T) \colon f(x) = h(x) \,\, \forall \ x \in T \setminus \interior{F}\}$. In terms of $\mathit{a}$, we have that $\mathit{a} = \sigma(\cup_{n \in \mathbb{N}}\mathit{a}_{n})$,
with 
\begin{gather*}
    \mathit{a}_{n} = \sigma\big(G_{F}^{n} \cup \{E\}\big), \quad G_{F}^{n} = \{\pi_{s_{1}, \dotsc, s_{n}}^{-1}(B_{n}) \colon B_{n} \in \mathfrak{B}(\mathbb{R}^{n})\}.
\end{gather*}
We write $\mathcal{U}_{n}$ for the collection of open balls in $\mathbb{R}^{n}$ with rational radii, centered around points in $\mathbb{Q}^{n}$.
Then, $\mathit{a} = \sigma(\Gamma_{E})$,
where 
\begin{gather*}
    \Gamma_{E} = \{\pi_{s_{1}, \dotsc, s_{n}}^{-1}(U_{n}) \colon U_{n} \in \mathcal{U}_{n}, n \in \mathbb{N}\} \cup \{E\}.
\end{gather*}
We remark that the subset $E$ is s.t.\ the sets $\pi_{s_{1}, \dotsc, s_{n}}^{-1}(U_{n})$, $U_{n} \in \mathcal{U}_{n}$, convey enough information to separate distinct elements of $E$ (Figure~\ref{fig:visualE}). Formally, we have the following lemma. 

\begin{figure}[ht]
  \centering
  \includegraphics[width=0.7\textwidth]{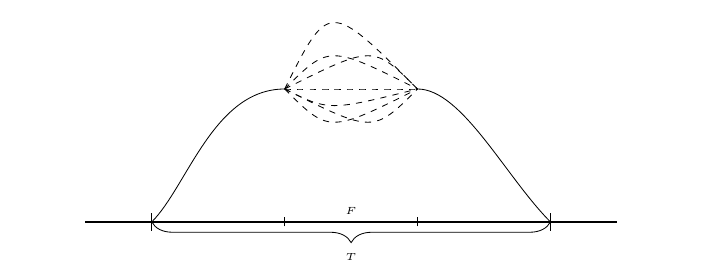}
  \caption{Schematic representation of $E$.}
  \label{fig:visualE}
\end{figure}

\begin{lemma} \label{lemma2}
    $\Gamma_{E}$ separates the points of $E$.
\end{lemma}
\begin{proof}
    Let $f,g \in E$ s.t.\ $f \neq g$. Then, since $f(v) = g(v)$ for any $v \in T \setminus \interior{F}$, there exists $p$ in $\interior{F}$ s.t.\ $f(p) \neq g(p)$. Since $f$, $g$ are continuous and $\overline{S} = F$, there exists $i \in \mathbb{N}$ s.t.\ $f(s_{i}) \neq g(s_{i})$. Then,
    \begin{gather*}
        \big(f(s_{1}), \dotsc, f(s_{i}), \dotsc, f(s_{n})\big) = v_{n} \neq w_{n} = \big(g(s_{1}), \dotsc, g(s_{i}), \dotsc, g(s_{n})\big).
    \end{gather*}
    Hence, there exists $U \in \mathcal{U}_{n}$ s.t.\ $v_{n} \in U_{n}$ but $w_{n} \notin U_{n}$. 
\end{proof}
Using the previous lemma, we readily see that $\Gamma_{E}$ is a basis for $(E, \mathit{a}, P_{\mathit{a}})$, where $P_{\mathit{a}}$ is the restriction of $P$ to $\mathit{a}$. Thus, since $E$ is clearly not countable, using Theorem~\ref{thm2}, \eqref{eq:martin} can not be satisfied.
\end{example}

For further illustration, the following example treats the case where the target space is a sequence space over a countable set, equipped with a product measure. 

\begin{example} \label{example4}
    Let $X = \mathbb{Z}^{\mathbb{N}} = \{x \colon x \colon \mathbb{N} \to \mathbb{Z}\}$ and $\pi_{1, \dotsc, n}(x) = (x(1), \dotsc, x(n))$, $x \in \mathbb{Z}^{\mathbb{N}}$, $n \in \mathbb{N}$. Let 
    \begin{gather*}
        \Sigma = \sigma\bigg(\bigcup_{n \in \mathbb{N}}\sigma\Big(\{\pi_{1, \dotsc, n}^{-1}(s_{1}, \dotsc,s_{n}) \colon s_{1}, \dotsc, s_{n} \in \mathbb{Z}\}\Big)\bigg),
    \end{gather*}
    and $P$ be a measure on $(\mathbb{Z}^{\mathbb{N}}, \Sigma)$ s.t.\ for any $s_{1}, \dotsc, s_{n} \in \mathbb{Z}$,
    \begin{gather*}
        P(\pi_{1, \dotsc, n}^{-1}(s_{1}, \dotsc,s_{n})) = p(\{s_{1}\}) \cdot \dotsc \cdot p(\{s_{n}\}),
    \end{gather*}
    where $p$ is a measure on $(\mathbb{Z}, 2^{\mathbb{Z}})$. The space $(\mathbb{Z}^{\mathbb{N}}, \Sigma,P)$ is a Lebesgue space with basis
    \begin{gather*}
        \big\{\pi_{1, \dotsc, n}^{-1}(s_{1}, \dotsc,s_{n}) \colon s_{1}, \dotsc, s_{n} \in \mathbb{Z}, n \in \mathbb{N}\big\}.
    \end{gather*}
    Compare also to Example~1.3 in \cite{entropy}. Let $S \subset \mathbb{Z}$ ($S$ contains at least two elements) and 
    \begin{gather*}
        \mathit{a}_{n} = \sigma\Big(\{\pi_{1, \dotsc, n}^{-1}(s_{1}, \dotsc,s_{n}) \colon s_{1}, \dotsc, s_{n} \in S\}\Big).
    \end{gather*}
    Then, \eqref{eq:martin} can not be satisfied with $\mathit{a} = \sigma(\cup_{n \in \mathbb{N}}\mathit{a}_{n})$. This follows from Theorem~\ref{thm2} with $E = S^{\mathbb{N}}$ and $\Gamma_{E} = \{\pi_{1, \dotsc, n}^{-1}(s_{1}, \dotsc,s_{n}) \colon s_{1}, \dotsc, s_{n} \in S, n \in \mathbb{N}\}$.
\end{example}

To conclude this section, the following example covers a case in which \eqref{eq:martin} is satisfied, $\mathit{a}$ is non-trivial, and the target space is uncountable.

\begin{example} \label{examplePro1}
    Consider $\{0,1\}^{\mathbb{N}} = \{x \colon x \colon \mathbb{N} \to \{0,1\}\}$ and let $\pi_{1, \dotsc, n}$ and $\Sigma$ be as in the previous example where $\mathbb{Z}$ is replaced with $\{0,1\}$.
    Let $y \in \{0,1\}^{\mathbb{N}}$ be s.t.\ $y(k) = 1$ for any $k \in \mathbb{N}$. Define
    \begin{gather*}
        P_{y}(A) = \delta_{y}(A) = \begin{cases}
            1, & \text{if $y \in A$,}\\
            0, & \text{otherwise,} 
        \end{cases} \quad A \in \Sigma.
    \end{gather*}
    Let $\mathit{a}_{n} = \sigma(\{x \colon x(k) = 1 \ \forall \, k \geq n\})$. We set $\mathit{a} = \sigma(\cup_{n \in \mathbb{N}}\mathit{a}_{n})$ and notice that $(\mathit{a}_{n})$ is increasing. Let $\xi \in Z_{0}^{\mathit{a}}$ and write $\xi = \{A_{y}\} \cup \{A \in \xi \colon y \notin A\}$, where $A_{y}$ is the unique atom of $\xi$ which contains $y$. Then, for any $A \in \xi$, by definition of $P_{y}$,
    \begin{gather*}
        P_{y}(A \mid \mathit{a}_{1}) = P_{y}(A \cap \{y\})\mathbbm{1}_{\{y\}}.
    \end{gather*}
    Therefore,
    \begin{gather*}
        E_{y}\big[-\sum_{A \in \xi}P_{y}(A\mid \mathit{a}_{1})\log P_{y}(A\mid \mathit{a}_{1}) \big] = E_{y}\big[-P_{y}(\{y\})\mathbbm{1}_{\{y\}}\log P_{y}(\{y\})\mathbbm{1}_{\{y\}}\big] = 0.
    \end{gather*}
\end{example}

\section{Maximal conditional entropy and families of martingales} \label{se:martingales}

\subsection{On the uniform convergence of families of martingales} \label{se:UniformConvergence1}

In this section, we assume that $\Theta$ is some nonempty set. The following two results arise as a consequence of \eqref{eq:martin}, they are proven in \cite{martin}.  

\begin{theorem}[Corollary of Theorem~1 in \cite{martin}] \label{corollary1}
    Let $(\mathit{a}_{n})$ be an increasing sequence of sub $\sigma$-fields s.t.\ $\mathit{a} = \sigma(\cup_{n \in \mathbb{N}}\mathit{a}_{n})$. For each $\theta \in \Theta$, let $(X_{n}^{\theta})$ be a martingale adapted to $(\mathit{a}_{n})$. Suppose that there exists $M \in \mathbb{R}$ s.t.\ for any $x \in X$, $\sup_{\{\theta \in \Theta, n \in \mathbb{N}\}}\lvert X_{n}^{\theta}(x) \rvert \leq M$. Then, if \eqref{eq:martin} is satisfied, there exists a random variable $X^{\theta}$ s.t.\
    \begin{gather*}
        \sup_{\theta \in \Theta}E\big[\lvert X_{n}^{\theta} - X^{\theta} \rvert\big] \xrightarrow[]{n \to \infty} 0.
    \end{gather*}   
\end{theorem}
 
\begin{theorem}[Corollary of Theorem~2 in \cite{martin}] \label{corollary2}
    Let $(\mathit{a}_{n})$ be an increasing sequence of sub $\sigma$-fields s.t.\ $\mathit{a} = \sigma(\cup_{n \in \mathbb{N}}\mathit{a}_{n})$. For each $\theta \in \Theta$, let $(X_{n}^{\theta})$ be a martingale adapted to $(\mathit{a}_{n})$. Assume that there exists $M \in \mathbb{R}$ s.t.\ for any $x \in X$, $\sup_{\{\theta \in \Theta, n \in \mathbb{N}\}}\lvert X_{n}^{\theta}(x) \rvert \leq M$. Then, if $\Theta$ is countable and \eqref{eq:martin} is satisfied, there exists a random variable $X^{\theta}$ s.t.\
    \begin{gather*}
        P\Big(\sup_{\theta \in \Theta}\lvert X_{n}^{\theta} - X^{\theta} \rvert \xrightarrow[]{n \to \infty} 0 \Big) = 1.
    \end{gather*}   
\end{theorem}

As a consequence of Theorem~\ref{corollary2}, we prove the following result, which allows for not necessarily countable $\Theta$.

\begin{corollary} \label{corollary3}
    Suppose that $\Theta \subset \mathbb{R}^{p}$ is closed. Let $(\mathit{a}_{n})$ be an increasing sequence of sub $\sigma$-fields s.t.\ $\mathit{a} = \sigma(\cup_{n \in \mathbb{N}}\mathit{a}_{n})$. For each $\theta \in \Theta$, let $(X_{n}^{\theta})$ be a martingale adapted to $(\mathit{a}_{n})$. Assume that
    \begin{enumerate}[(i)]
        \item \label{item1:C3} there exists $M \in \mathbb{R}$ s.t.\ for any $x \in X$, $\sup_{\{\theta \in \Theta, n \in \mathbb{N}\}}\lvert X_{n}^{\theta}(x) \rvert \leq M$;
        \item \label{item2:C3} for any $n \in \mathbb{N}$ and $x \in X$, $\theta \mapsto X_{n}^{\theta}(x)$ is continuous;
        \item \label{item3:C3} Given any $\theta \in \Theta$, there exists a random variable $X^{\theta}$ which is s.t.\ for any $x \in X$, $\theta \mapsto X^{\theta}(x)$ is continuous and $P(X_{n}^{\theta} \xrightarrow[]{n \to \infty} X^{\theta}) = 1$.
    \end{enumerate}
    Then, if \eqref{eq:martin} is satisfied, it follows that $P(\sup_{\theta \in \Theta}\lvert X_{n}^{\theta} - X^{\theta} \rvert \xrightarrow[]{n \to \infty} 0 ) = 1$.
\end{corollary}

\begin{proof}
    Let $\Theta^{\prime} = \Theta \cap \mathbb{Q}^{p}$. Since \eqref{eq:martin} is satisfied, using (\ref{item1:C3}), it follows from Theorem~\ref{corollary2} that there exists a random variable $\widetilde{X}^{q}$, $q \in \Theta^{\prime}$, s.t.\ $P(\sup_{q \in \Theta^{\prime}}\lvert X_{n}^{q} - \widetilde{X}^{q} \rvert \xrightarrow[]{n \to \infty} 0) = 1$.
    Let 
    \begin{gather*}
        A = \bigcap_{q \in \Theta^{\prime}}\{x \colon X_{n}^{q}(x) \xrightarrow[]{n \to \infty} X^{q}(x)\}.
    \end{gather*}
    Using~(\ref{item2:C3}) and~(\ref{item3:C3}), we have that
    \begin{gather*}
        \big\{x \colon \sup_{q \in \Theta^{\prime}}\lvert X_{n}^{q}(x) - \widetilde{X}^{q}(x) \rvert \xrightarrow[]{n \to \infty} 0\big\} \cap A \subset \big\{x \colon \sup_{\theta \in \Theta}\lvert X_{n}^{\theta}(x) - X^{\theta}(x) \rvert \xrightarrow[]{n \to \infty} 0\big\},
    \end{gather*}
    which completes the proof.
\end{proof}

\subsection{On maximal conditional entropy and random fields} \label{se:UniformConvergence2}

This final section covers a specific case where \eqref{eq:martin} does not apply. However, the uniform convergence of a specific type of martingale plays a prominent role in obtaining consistent parameter estimators. Let $X = \{x \colon x \colon E \to \mathbb{R}\}$, where $E \subset \mathbb{R}^{d}$. We equip $X$ with the coordinate $\sigma$-field specified along coordinates which are dense in $E$. That is, we let $\mathcal{B} = \sigma(\cup_{n \in \mathbb{N}}\mathit{a}_{n})$, where 
\begin{gather*}
    \mathit{a}_{n} = \sigma\big(\{\pi_{s_{1}, \dotsc, s_{n}}^{-1}(B_{n}) \colon B_{n} \in \mathfrak{B}(\mathbb{R}^{n})\}\big), \quad \pi_{s_{1}, \dotsc, s_{n}}(x) = (x(s_{1}), \dotsc, x(s_{n})),
\end{gather*}
and $\{s_{i} \colon i \in \mathbb{N}\}$ is dense in $E$. On a probability space $(\Omega, \mathcal{F}, \mathbb{P})$, let $Y$ to be a random function with values in $X$ and distribution $P$ on $\mathcal{B}$. We assume that $(X, \mathcal{B}, P)$ is a Lebesgue space. Specific cases are $X = C(T)$ in which the coordinate $\sigma$-field actually signifies the Borel $\sigma$-field on $C(T)$ (cf.\ Example~\ref{example1}) or $X$ is a certain sequence space with $E = S = \mathbb{N}$ and $x \in X$ takes values in a countable subset of $\mathbb{R}$ (cf.\ Example~\ref{example4}). Let $Y_{n} = (Y_{s_{1}}, \dotsc, Y_{s_{n}})$. Then, 
\begin{equation} \label{eq:RandomVector}
    \sigma(Y_{n}) = \{Y_{n}^{-1}(B_{n}) \colon B_{n} \in \mathfrak{B}(\mathbb{R}^{n})\} = \{Y^{-1}(A) \colon A \in \mathit{a}_{n}\}.
\end{equation}
Suppose that $\{P_{\theta} \colon \theta \in \Theta\}$ is a family of distributions for $P = P_{\theta_{0}}$, $\theta_{0} \in \Theta$, where $\Theta = \prod_{j=1}^{p}[a_{j}, b_{j}] \subset \mathbb{R}^{p}$, $a_{j} < b_{j}$, $j = 1, \dotsc,p$. Given $\theta \in \Theta$, let $P_{\theta}^{n}$ be the restriction of $P_{\theta}$ to $\mathit{a}_{n}$. We assume that for any $n \in \mathbb{N}$ and for any $\theta \in \Theta$, $P_{\theta}^{n} \ll P_{\theta_{0}}^{n}$ on $\mathit{a}_{n}$ and write
\begin{gather*}
    \varphi_{n}^{\theta}(x) = \frac{P_{\theta}^{n}(dx)}{P_{\theta_{0}}^{n}(dx)}, \quad x \in X,  
\end{gather*}
for the Radon–Nikodym derivative of the measure $P_{\theta}^{n}$ with respect to $P_{\theta_{0}}^{n}$. At this point, two remarks are helpful. First, if we let $p_{\theta}^{n}$, $\theta \in \Theta$, be the respective family of distributions of the random vector $Y_{n}$ on $\mathfrak{B}(\mathbb{R}^{n})$, then because of \eqref{eq:RandomVector}, the absolute continuity of the measures $P_{\theta}^{n}$ and $P_{\theta_{0}}^{n}$ on $\mathit{a}_{n}$ can be deduced from the absolute continuity of the measures $p_{\theta}^{n}$ and $p_{\theta_{0}}^{n}$ on $\mathfrak{B}(\mathbb{R}^{n})$. Second, if we define
\begin{gather*}
    L_{n}^{\theta}(x) = \phi_{n}^{\theta}\big(\pi_{s_{1}, \dotsc, s_{n}}(x)\big), \quad x \in X,
\end{gather*}
where $\phi_{n}^{\theta}(y) = p_{\theta}^{n}(dy)/p_{\theta_{0}}^{n}(dy)$, $y \in \mathbb{R}^{n}$, then again by \eqref{eq:RandomVector}, we have that
\begin{gather*}
    P_{\theta}^{n}(A) = \int_{A}L_{n}^{\theta}(x)P_{\theta_{0}}^{n}(dx), \quad A \in \mathit{a}_{n}.
\end{gather*}
Thus, since we assume that $P_{\theta}^{n} \ll P_{\theta_{0}}^{n}$ on $\mathit{a}_{n}$ with derivative $\varphi_{n}^{\theta}$, we have that $L_{n}^{\theta} = \varphi_{n}^{\theta}$ with $P_{\theta_{0}}^{n}$ probability one on $X$. Notice that by definition, $L_{n}^{\theta}$ is $\mathit{a}_{n}$ measurable for any $n \in \mathbb{N}$. Hence, $P_{\theta_{0}}(L_{n}^{\theta} = \varphi_{n}^{\theta})=1$. Using Theorem~\ref{thm1}, if $(X, \mathcal{B}, P_{\theta_{0}})$ is a Lebesgue space, we know that \eqref{eq:martin} is not satisfied for this particular case where $\mathit{a} = \mathcal{B}$. It is shown in Section~1 of chapter VII in \cite{StochasticProcesses} that for any $\theta \in \Theta$, $(\varphi_{n}^{\theta})$ is a martingale on $(X, \mathcal{B}, P_{\theta_{0}})$, adapted to the filtration $(\mathit{a}_{n})$. In particular, since $L_{n}^{\theta}$ is $\mathit{a}_{n}$ measurable for any $n \in \mathbb{N}$ and $P_{\theta_{0}}(L_{n}^{\theta} = \varphi_{n}^{\theta})=1$, it follows that for any $\theta \in \Theta$, $(L_{n}^{\theta})$ is a martingale on $(X, \mathcal{B}, P_{\theta_{0}})$, adapted to the filtration $(\mathit{a}_{n})$. The following proposition is known, it is based on an argument proposed by \cite{wald} (cf.\ Theorem~2 in \cite{wald}).

\begin{proposition} \label{proposition2}
    Assume that for any $n \in \mathbb{N}$ and $\theta \in \Theta$, $p_{\theta}^{n}$ is absolutely continuous with respect to the Lebesgue measure on $\mathfrak{B}(\mathbb{R}^{n})$ with probability density function $f_{n}^{\theta}(y)$, $y \in \mathbb{R}^{n}$. Let $\hat{\theta}_{n} \colon \Omega \mapsto \Theta$, $n \in \mathbb{N}$, be a sequence of random vectors on $(\Omega, \mathcal{F}, \mathbb{P})$. Suppose that
    \begin{enumerate}[(i)]
        \item \label{item1:C5} for any $\theta \in \Theta$ and $n \in \mathbb{N}$, 
        \begin{equation} \label{eq:consistency1}
             f_{n}^{\hat{\theta}_{n}(\omega)}(Y_{n}(\omega)) \geq f_{n}^{\theta}(Y_{n}(\omega)), \quad \omega \in \Omega;
        \end{equation}
        \item \label{item2:C5} for any $\varepsilon > 0$ s.t.\ $B(\theta_{0}, \varepsilon) \subset \Theta$,
        \begin{equation}\label{eq:consistency5}
             \mathbb{P}\bigg(\sup\bigg\{\frac{f_{n}^{\theta}(Y_{n})}{f_{n}^{\theta_{0}}(Y_{n})} \colon \theta \in \Theta \setminus B(\theta_{0}, \varepsilon)\bigg\} \xrightarrow[]{n \to \infty} 0\bigg) = 1.
        \end{equation}
    \end{enumerate}
    Then,
    \begin{equation}\label{eq:consistency4}
      \mathbb{P}\Big(\hat{\theta}_{n} \xrightarrow[]{n \to \infty} \theta_{0}\Big) = 1.
    \end{equation}
\end{proposition}

\begin{remark} \label{remark1}
Suppose that $p_{\theta}^{n}$ is absolutely continuous with respect to the Lebesgue measure on $\mathfrak{B}(\mathbb{R}^{n})$ with probability density function $f_{n}^{\theta}(y)$. Then, \eqref{eq:consistency5} of Proposition~\ref{proposition2} is satisfied if for any $\varepsilon > 0$ s.t.\ $B(\theta_{0}, \varepsilon) \subset \Theta$,  
\begin{gather*}
    P_{\theta_{0}}\big(\sup\{L_{n}^{\theta} \colon \theta \in \Theta \setminus B(\theta_{0}, \varepsilon)\} \xrightarrow[]{n \to \infty} 0\big) = 1.
\end{gather*}
An estimator $(\hat{\theta}_{n})$ which satisfies \eqref{eq:consistency1}, is referred to as a sequence of maximum likelihood estimators for $\theta_{0}$. If $(\hat{\theta}_{n})$ satisfies \eqref{eq:consistency4}, $(\hat{\theta}_{n})$ is said to be (strongly) consistent for $\theta_{0}$. However, since \eqref{eq:martin} is not satisfied, Corollary~\ref{corollary3} does not apply. We note that \eqref{eq:consistency5} can be recovered using independence (cf.\ \cite{wald}) or impose further regularity conditions on $\theta \mapsto f_{n}^{\theta}(y)$ (see for instance the proof of Theorem~3 in \cite{zhang} or more lately Theorem~4.1 in \cite{hediger}).   
\end{remark}



\section*{Acknowledgments}
I thank my supervisor Reinhard Furrer for all the helpful comments, questions and thoughts. This study was supported by the Swiss National Science Foundation SNSF-175529.

\bibliographystyle{amsplain} 
\bibliography{references.bib}

\end{document}